\documentclass[cmbright]{amsart}
\usepackage{moreverb}

\newcommand\BibTeX{{\rmfamily B\kern-.05em \textsc{i\kern-.025em b}\kern-.08em
T\kern-.1667em\lower.7ex\hbox{E}\kern-.125emX}}

\usepackage{amsmath,amssymb, amsfonts}
\usepackage{graphicx}
\usepackage{caption, subcaption}
\usepackage{pdfsync}
\usepackage{color,bm}
\usepackage{acronym}
\usepackage{fullpage}

\newcommand{\la}{\lambda}
\newcommand{\Om}{\varOmega}

\newcommand{\E}{\mathcal{E}}
\newcommand{\ii}{\mathrm{i}}
\newcommand{\pa}{\bar{\partial}}

\newcommand{\V}{\mathbb{V}}

\newcommand{\Nb}{\mathbb{N}}

\newcommand{\T}{\mathcal{T}}
\newcommand{\Pro}{\mathcal{P}}

\acrodef{CNLS}{cubic nonlinear Schr\"odinger}

\makeatletter
\newcommand{\sech}{\mathop{\operator@font sech}}
\newcommand{\sign}{\mathop{\operator@font sign}}
\makeatother

\begin{document}

\title{On the reflection of solitons of  the cubic nonlinear Schr\"{o}dinger equation\footnotemark[2]}

\author{Theodoros Katsaounis}
\address[Theodoros Katsaounis]{King Abdullah University of Science and Technology(KAUST)\\
Thuwal, Kingdom of Saudi Arabia \\
 \& IACM--FORTH, Heraklion, GREECE \\
 \& Dept. of Math \& Appl. Mathematics, Univ. of Crete, Greece
}
\email{theodoros.katsaounis@kaust.edu.sa}
\author{Dimitrios Mitsotakis}
\address[Dimitrios Mitsotakis]{Victoria University of Wellington, School of Mathematics and Statistics, New Zealand \\
 \& IACM--FORTH, Heraklion, GREECE \\
 \& Dept. of Math \& Appl. Mathematics, Univ. of Crete, Greece
}
\email{dimitrios.mitsotakis@vuw.ac.nz}

%
%

\begin{abstract}
In this paper we perform a numerical study on the interesting phenomenon of soliton reflection of solid walls. We consider the 2D cubic nonlinear Schr\"odinger equation as the underlying mathematical model and we use  an implicit-explicit type Crank-Nicolson finite element scheme for its numerical solution. After verifying the perfect reflection of the solitons on a vertical wall, we present the imperfect reflection of a dark soliton on a diagonal wall. 
\end{abstract}

\maketitle


\vspace{-6pt}

\section{Introduction}
In this paper we study numerically  the phenomenon of reflection of bright and dark solitons of walls. To this effect, 
we consider the initial value problem of the \acf{CNLS}  equation 
%
\begin{equation}
\label{NLS} \left \{
\begin{aligned}
& \ii u_t +\Delta u + \la |u|^2 u=0    &&\quad\mbox{in ${\Omega}\times (0,T]$,}& \\
&u(\cdot,0)=u_0 &&\quad\mbox{in ${\Omega}$},&
\end{aligned}\right.
\end{equation}
where we assume that  $\Om\subset\mathbb{R}^2$, is a bounded, convex, polygonal domain,  $T<\infty$, $\la\in\mathbb{R}$.  For $\lambda \le0$, problem \eqref{NLS} is known as the \emph{defocusing cubic nonlinear Schr\"odinger  equation}, while If $\lambda>0$,  problem \eqref{NLS} is called the \emph{focusing cubic nonlinear Schr\"odinger equation}.
The \acs{CNLS} equation \eqref{NLS} is used as a mathematical model in various applications such as  nonlinear optics and lasers,  water waves, quantum hydrodynamics  and Bose-Einstein condensates,  \cite{Sulem}. The nature of soliton reflection to walls requires that \eqref{NLS} is augmented with boundary conditions. In our study we consider three type of boundary conditions, a) \emph{zero Dirichlet} \eqref{DBC}  or  b) \emph{zero Neumann} \eqref{NBC} condition, in the whole boundary $\partial \Omega$ of the domain, or c) \emph{zero Dirichlet} and \emph{zero Neumann} on disjoint parts of the boundary \eqref{MBC}, 
\begin{align}
u &= 0 \ \text{ on } \ \partial \Omega, \label{DBC} \\
\frac{\partial u}{\partial n} & = 0 \ \text{ on }\  \partial \Omega, \label{NBC} \\
u & = 0 \ \text{ on }\  \partial \Omega^D \ \text{ and } \ \frac{\partial u}{\partial n} = 0 \ \text{ on } \partial \Omega^N,  \quad     \partial \Omega =  \partial \Omega^D \cup \partial \Omega^N, \quad \partial \Omega^D \cap \partial \Omega^N = \emptyset, \label{MBC}
\end{align}
%
A standard calculation, c.f. \cite{C},  shows that the initial value problem  \eqref{NLS}, augmented with any of the aforementioned boundary conditions \eqref{DBC}, \eqref{NBC} or \eqref{MBC}, conserves two physical  quantities: mass $M(\cdot)$ and energy $E(\cdot)$.   In particular, we have 
%
\begin{equation}
\label{CL}
M(t):= \|u(t)\|^2, \quad  E(t):= \frac{1}{2}\|\nabla u(t)\|^2 - \frac{\la}{4}\|u(t)\|^{4}_{L^{4}}, \ \text{ then }\ M(t) = M(0),  \ \text{ and } \  E(t) = E(0) \quad t\ge 0,
\end{equation}
where $\|\cdot\|$,  $\|\cdot\|_{L^q}$ denote the $L^2$ and $L^q-$norms in $\Om$, respectively.

In this numerical study we focus on two aspects: a) we evaluate a  Crank-Nicolson relaxation method for 2D domains discretized by completely unstructured grids and b) we study various reflections of solitons on walls for \acs{CNLS} type \eqref{NLS} of equations either focusing or defocusing.

Nonlinear Schr\"odinger type equations can describe waves in optical fibers, \cite{Argawal}, as well as rogue waves in the ocean, cf. e.g. \cite{Osborne}. In both cases the study of the interaction of the waves with structures imposes the study of the reflection of solitons on walls.
On the other hand, it is well known that in one space dimension \acs{CNLS} equation \eqref{NLS} is an integrable system, thus one can compute analytically exact solutions describing soliton reflections. However,  in two space dimensions system \eqref{NLS} is not integrable thus one has to rely on numerical methods for computing and studying such reflective phenomena.

The numerical method used here  is based on the standard  finite element method for the spatial discretization and the \emph{relaxation Crank-Nicolson} scheme as a time marching mechanism.  The relaxation Crank-Nicolson scheme was introduced by Besse \cite{Besse} can be viewed as a linearization method for the \acs{CNLS} and as such it avoids the computationally expensive solution of a nonlinear equation  at each time step of the algorithm. Moreover, the relaxation scheme exhibits mass conservation, same like the standard Crank-Nicolson scheme, thus preserving the mass conservation property of the continuous problem, cf. \eqref{CL}. However,  the relaxation scheme does not preserve the energy $E(\cdot)$ in two space dimensions. 
The relaxation Crank-Nicolson method is also used  in \cite{KK} where optimal a posteriori error estimates for models of type \eqref{NLS} are obtained. Based on these a posteriori error estimates one can derive a space-time adaptive algorithm which will be able to capture all interested features of the solution with substantial reduction of  the overall  computational cost compared to that of uniform grids. A space-time adaptive algorithm based on a posteriori error estimates was developed in \cite{KK0} for  the linear Schr\"{o}dinger equation.  
The numerical experiments reported in Section 3 are good examples where such a space-time adaptive algorithm will be very beneficial. 

Due to the integrability properties of the \acs{CNLS} equation, the reflection of a soliton can be studied analytically only in 1D, \cite{Biondini1,Biondini2,Tarasov}. In this paper after verifying the order of accuracy of the numerical method in space and time, we validate the efficiency of the numerical method by studying first the perfect (elastic) reflection of dark and bright solitons on vertical walls using either \eqref{NBC} or \eqref{MBC} as boundary conditions, \cite{Biondini1}. Finally, we show that the reflection of a dark soliton is not perfect (inelastic) when the soliton collides on the wall at an angle.

This paper is organized as follows: in Section 2 we describe briefly the relaxation Crank-Nicolson finite element method and state its approximation properties. Section 3 contains numerical results which a) validate the numerical method and b) study the phenomenon of soliton reflection.  

\section{The numerical Method}
We describe now briefly the numerical method used in this study.  
We consider a uniform partition   $t_n=n k $ of $[0,T]$ where $k = T/N$ and $I_n:=(t_n,t_{n+1}]$, $0\le n\le N-1,$ denote the fixed time step and subintervals of $[0,T],$ respectively. 
For the spatial discretization, we consider a family of conforming, shape regular triangulations $\{\T\}$ of $\Om$. 
For an element $K\in\T$, we denote  by $h_K$ its diameter and let $h=\max_K h_K$. We also let  $\mathbb{P}^r$ denotes the space of polynomials in two variables of degree at most $r$. Then to the  triangulation $\T$ we associate a finite element space $\V_r$ which definition will depend upon the type of boundary condition we choose to work with.  
In particular we define
$$\V_r:=\{\chi\in \V(\Om):\forall K\in\T,\, \chi|_K\in\mathbb{P}^r\},$$
where $\V=H_0^1(\Om)$ for homogeneous Dirichlet condition \eqref{DBC}, $\V=H^1(\Om)$ for the homogeneous Neumann boundary condition \eqref{NBC} and $\V=H_0^1(\Om^D)$ in the case of the mixed type boundary condition \eqref{MBC}.


We can define now the  relaxation Crank-Nicolson-Galerkin-type fully discrete scheme. Let  $\pa U^n:=(U^{n+1}-U^n)/k $ and $U^{n+\frac12}:=(U^{n+1}+U^n)/2$. Then we seek approximations $U^n\in\V_r$ to $u(t_n)$ such that, for $0\le n\le N-1$,
\begin{equation}
\label{CNrelaxfull} \left \{
\begin{aligned}
& \langle \frac12(\Phi^{n+\frac12}+\Phi^{n-\frac12}),\chi\rangle =\langle |U^n|^{2},\chi\rangle,   \quad\forall \chi\in \V_r,  \\
&\ii \langle \pa U^n, \chi\rangle - \langle\nabla U^{n+\frac12},  \nabla\chi\rangle + \la \langle \Phi^{n+\frac12} U^{n+\frac12}, \chi \rangle = 0,  \quad\forall \chi\in \V_r , 
\end{aligned}
\right.
\end{equation}
where $\langle\cdot,\cdot\rangle$  denotes the $L^2-$ inner product,  $\Phi^{-\frac12}=\Pro\left(|u_0|^{2}\right)$, $U^0=\Pro u_0$ with $\Pro$ being the $L^2-$projection $\Pro:L^2\to\V_r$.   Formally we expect the method to be second order accurate in time, and of $r+1$-order accurate in space, which can be expressed by an a priori error estimate of the following form 
\begin{equation}
\label{error}
\max_{0\le n\le N} \|U^n - u(t_n)\| \le C\left(h^{r+1} + k^2\right),
\end{equation}
where $C$ is a constant depending on the exact solution $u$ of \eqref{NLS} and data of the problem, but it is independent of $h$ and $k$. A similar estimate was proven rigorously in \cite{KAD} and in the form of a posteriori error bound in \cite{KK}. In \cite{Besse} an analogous a priori error estimate was obtained using finite differences. 

At each time step $t_n$, given an approximation $U^n\in\V_r$ and $\Phi^{n-\frac12}$, the algorithm proceeds first by computing the new value  $\Phi^{n+\frac12}$. The computational cost of this update is relatively low since it involves only the cost of the projection of the right hand side in the first relation of \eqref{CNrelaxfull}.  The cost of such projection amounts to  the solution of a real linear system with the mass matrix and real right hand side since $\Phi^{n+\frac12}$ is real valued. The second part of the algorithm \eqref{CNrelaxfull}, involves the update of the right hand side via a projection and the solution of a complex-valued banded linear system. An alternative approach for solving the linear systems is by using an appropriate conjugate gradient type method, since the system matrices in both steps are positive definite. Therefore the algorithmic complexity of method \eqref{CNrelaxfull} at each time step is at most of $O(m^2)$ where $m$ denotes the size of the matrix, thus making the method computationally very attractive.

\section{Numerical Experiments}
In this section we first validate numerically the method \eqref{CNrelaxfull} by means of verifying the formal convergence rates in space and time of the formal error estimate \eqref{error}. Furthermore, we present some results concerning the head on or oblique soliton reflection of solid walls.

In all, expect one, of  the numerical experiments, we used unstructured triangulations of Delaunay type produced by the Bowyer-Watson and Chew's second algorithm, \cite{Shewchuk}.  Delaunay triangulations have several advantages compared to uniform triangulations with the first being that are better suited for general geometries of the domain. 
Further, one of  characteristics  of the Delaunay triangulations is that maximizes the minimum interior angle of the triangles and on the same time minimizes the maximum interior angle of the triangles, which guarantees triangles of good quality. This characteristic is preserved even when the triangulation is refined and/or coarsened locally and the resulting triangulation is also Delaunay. However this is not the case with uniform triangulations since local refinement and/or coarsening can deteriorate the quality of triangles rather fast and produce triangles which are elongated and "thin" with very small angles.   Uniform triangulations may favour specific geometric directions in the domain which can influence the computed solution; such an issue is not present in an unstructured Delaunay triangulation.

\subsection{Method validation}
To validate the method numerically we perform a series of numerical experiments verifying the order of convergence as presented by the formal {\em a priori} error estimate \eqref{error}. To facilitate the process and be able to compute the exact error between the true and approximate solution, we choose to work with a solution of a non-homogeneous version of equation \eqref{NLS} and in particular we choose 
\begin{equation}
\label{exsol}
u(x,y,t) =  e^t x(1-x)y(0.5-y),
\end{equation}
which is an exact solution of \eqref{NLS} with appropriate right-hand side, with $\lambda=-2$, zero Dirichlet boundary conditions \eqref{DBC}, and $u_0(x,y)=u(x,t,0)$. For the computation of the errors in time we use quadratic finite elements ($r=2$) in space, the computational domain is the rectangle $\Om=[0,1]\times[0,0.5]$ and final time $T=4$. The domain is covered by an unstructured in general triangulation $\T$ of good quality.  To verify the error convergence rates we compute the \emph{experimental order of convergence(EOC)}. Let $\ell\in\Nb$ count the different realizations(runs) and let  $h_{\ell}, \ k_{\ell},\  \E_{\ell}$ be the spatial mesh size, time step and error respectively.  Choosing a very small spatial mesh size, we used an unstructured triangulation consisting of 192802 triangles, the corresponding spatial component of the error is negligible and the temporal EOC is computed as $\displaystyle \text{ EOC }=\frac{\log{(\E_{\ell+1} / \E_{\ell})}}{\log{(k_{\ell+1} / k_{\ell})}}$. The temporal EOC is found to be 2 and is presented Table \ref{temprate}. 

\begin{table}[htbp]
\vskip 20pt
\begin{center}
\begin{tabular}{cccccc}
\hline
$\ell$ & $1$ & $2$ & $3$ & $4$ & $5$\\
\hline
$k_{\ell}$ &  $0.5$ & $0.25$ & $0.08$ & $0.0625$ & $0.03125$ \\
$E_{\ell}$ & $9.02\times 10^{-3}$ & $2.27\times 10^{-3}$ & $2.35\times 10^{-4}$ & $1.43\times 10^{-4}$ & $3.49\times 10^{-5}$\\
EOC & -- & $1.9880$ & $1.9905$ & $2.0020$ & $2.0412$\\
\hline
\end{tabular}
\end{center}
\caption{Temporal numerical errors and experimental orders of convergence}\label{temprate}
\end{table}

\begin{figure}[htbp]
\vspace{-1.5cm}
\centering
\includegraphics[scale=0.4]{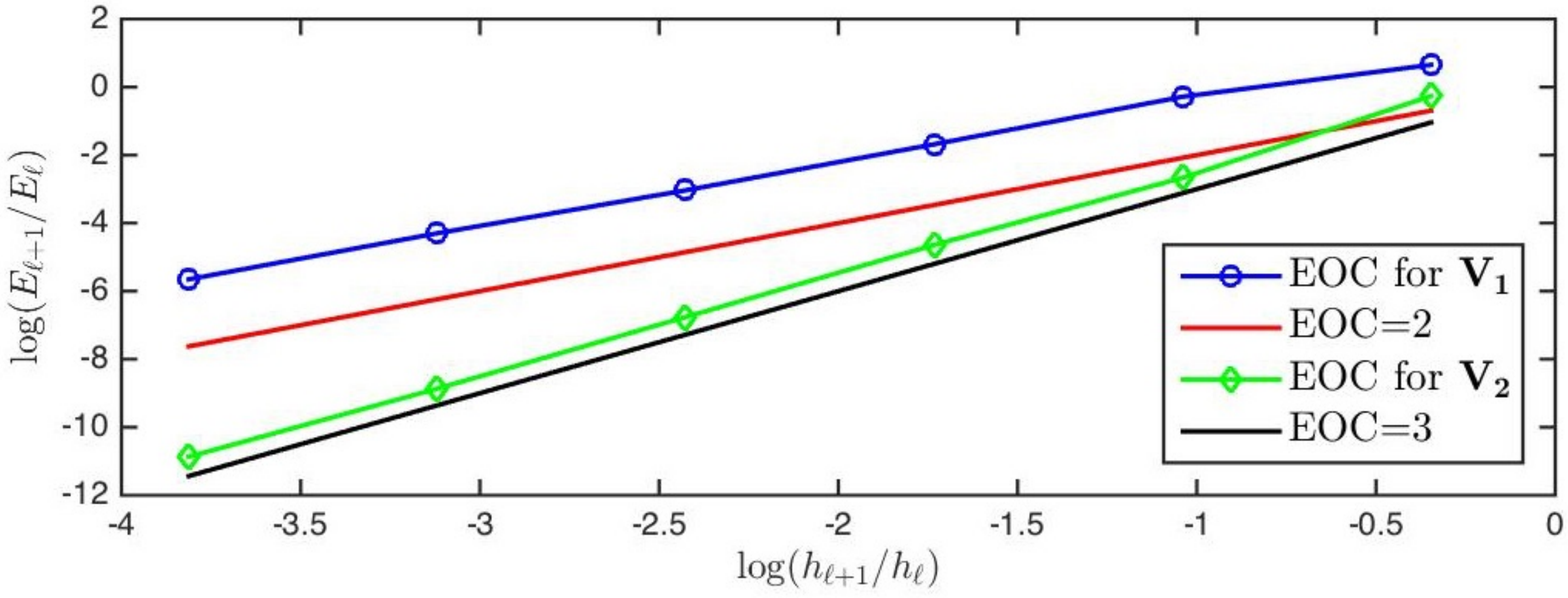}
\vspace{-1.75cm}
  \caption{Experimental orders of convergence for linear and quadratic elements.}
  \label{spatrate}
\end{figure}

The spatial errors can be computed in a similar manner. Specifically, in order to estimate the spatial convergence rates we take the domain to be $\Om=[0,2]\times [0,2]$ and exact solution $u(x,y,t)=e^t(1-\cos(2\pi x))\sin(2\pi y)$, with $u=0$ on the boundary,  along with the appropriate non-homogeneous term.
Now we take $T=0.1$ and $k=2\cdot 10^{-5}$, thus the temporal component of the error is negligible. The triangulation now is structured(uniform) and consisted of equal right-angle triangles with perpendicular sides of length $h=\sqrt{2/N}$, where $N$ is the number of triangles. The number of triangles we tested  were $N=32, 128, 512, 2048, 8192, 32768$.  The numerical experiments confirmed the expected orders of convergence. The numerical results are depicted in Figure \ref{spatrate}, using logarithmic scale in both axes. 

\subsection{Perfect Reflection of Solitons}

The focusing \acs{CNLS} equation ($\lambda=2$) has bright soliton solutions of the form $u_b(x,y,t)=\eta\sech\left[\eta(x+2\xi t)\right] e^{-i\theta}$ where $\theta=\xi x+(\xi^2-\eta^2)t$, while the defocusing \acs{CNLS} equation ($\lambda=-2$) admits dark soliton solutions of the form $u_d(x,y,t)=\eta\left[\cos\xi+{\rm i} \sin\xi\tanh(\sin\xi \eta(-x+2\eta\cos\xi t))\right]e^{-2i\eta^2 t}$, cf. e.g. \cite{Ablowitz}. The parameters $\xi$ and $\eta$ are chosen appropriately. 

In this section we verify that the reflection of a soliton of the \acs{CNLS} equation is \emph{perfect} when the soliton collides with a vertical wall at zero angle. The reflection is called \emph{perfect} if the reflected wave has the same shape as the original soliton but different direction of propagation. This behaviour has been studied analytically in \cite{Biondini1,Biondini2,Tarasov} for the integrable \acs{CNLS} equation in one space dimension. We performed two numerical tests using boundary conditions \eqref{NBC} and \eqref{MBC}.  Both type of boundary conditions will give perfect reflections but the interaction of the soliton with the boundary in different.  For the differences between the two reflections we refer to \cite{Biondini1,Biondini2}.

\begin{figure}[ht!]
\centering
\hspace{-2.5cm}
\begin{subfigure}[b]{0.55\textwidth}
\includegraphics[scale=0.35]{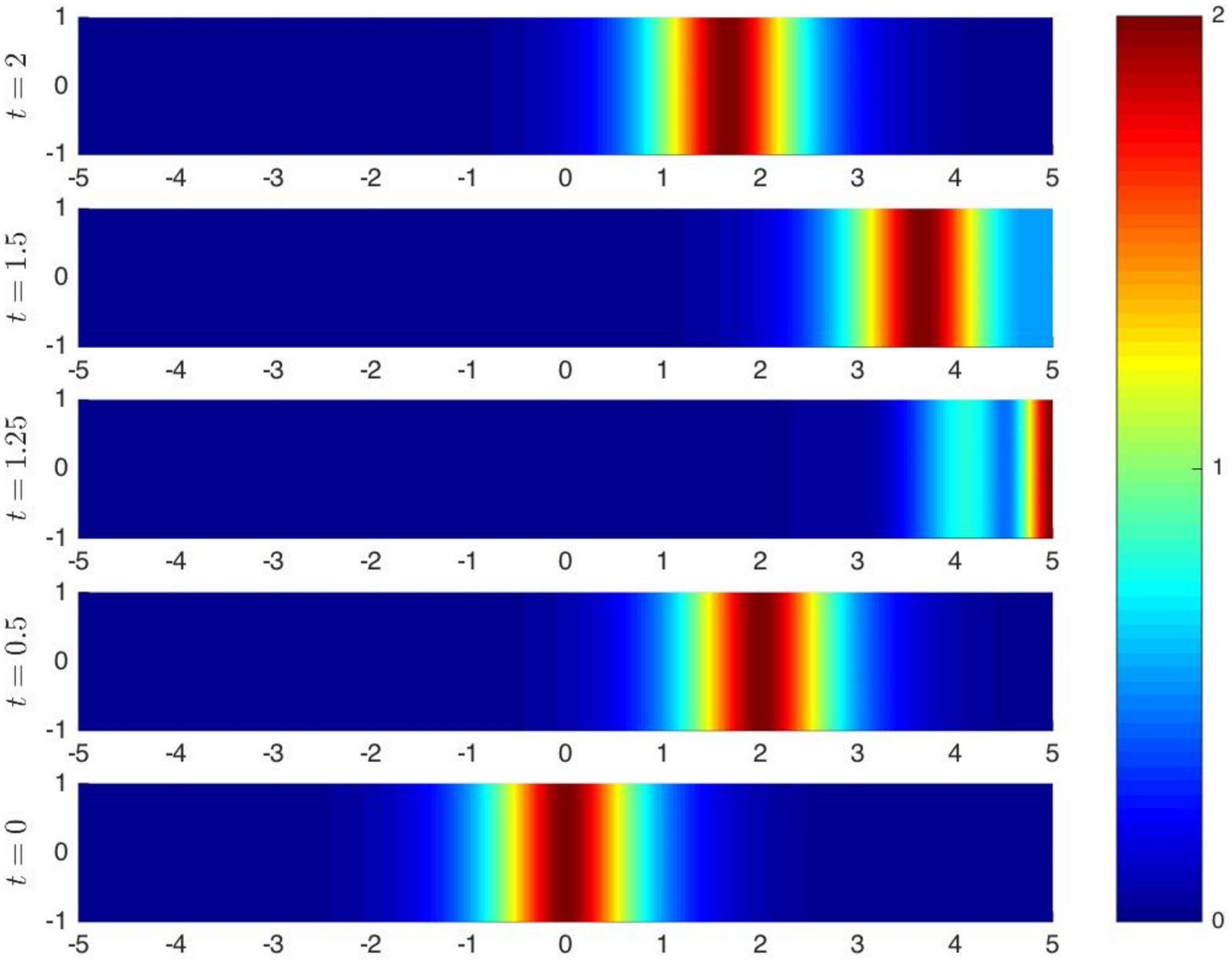}
\caption{Homogeneous Neumann boundary conditions \eqref{NBC}.}
\label{refl1a}
\end{subfigure}
\begin{subfigure}[b]{0.45\textwidth}
\includegraphics[scale=0.235]{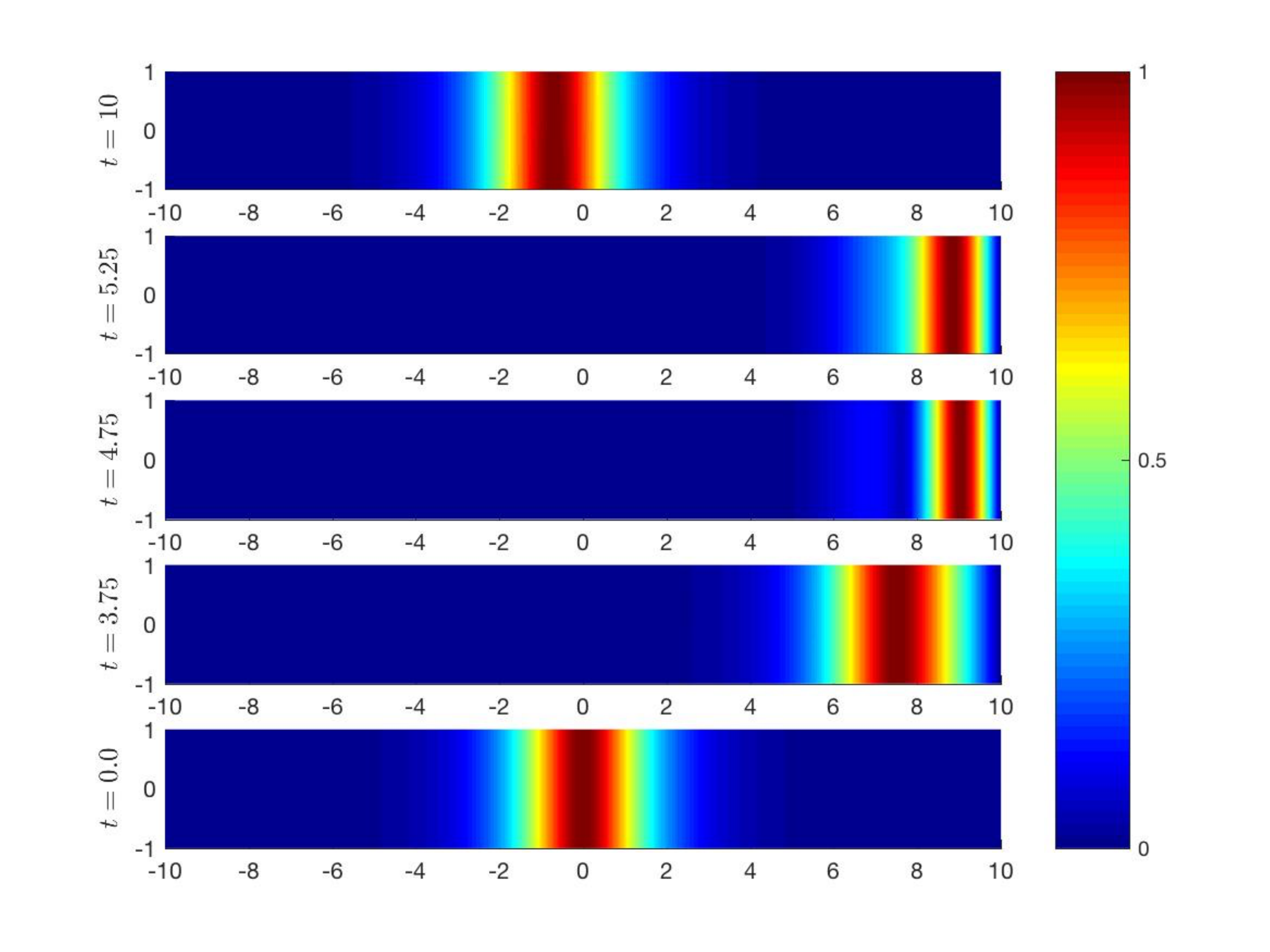}
\caption{Mixed boundary conditions \eqref{MBC}.}
\label{refl1b}
\end{subfigure} 
\caption{Perfect reflection of bright solitons}
\label{refl1}
\end{figure}

In the first test we study the reflection of a bright soliton for the focusing \acs{CNLS} equation with $\lambda=2, \ \eta=2$, $\xi=2$ and zero Neumann boundary conditions \eqref{NBC} in the domain $\Omega=[-5,5]\times[-1,1]$, while quadratic elements  were used.   The results of the perfect reflection are presented in Figure \ref{refl1a}, depicting the amplitude of the wave. It is known that the solitons of the \acs{CNLS} equation suffer by an instability of focusing type. In order to ensure that the propagation of the soliton remains stable during the simulation we took $74496$ triangles ensuring a very fine spatial grid and a small time-step $k=5\times 10^{-3}$. During the experiment the mass $M$ was conserved with value $7.9999999$ while the energy $E$ was conserved to $10.5$ up to $T=3$.

\begin{figure}[ht!]
\centering
\hspace{-2.5cm}
\begin{subfigure}[b]{0.55\textwidth}
\includegraphics[scale=0.33]{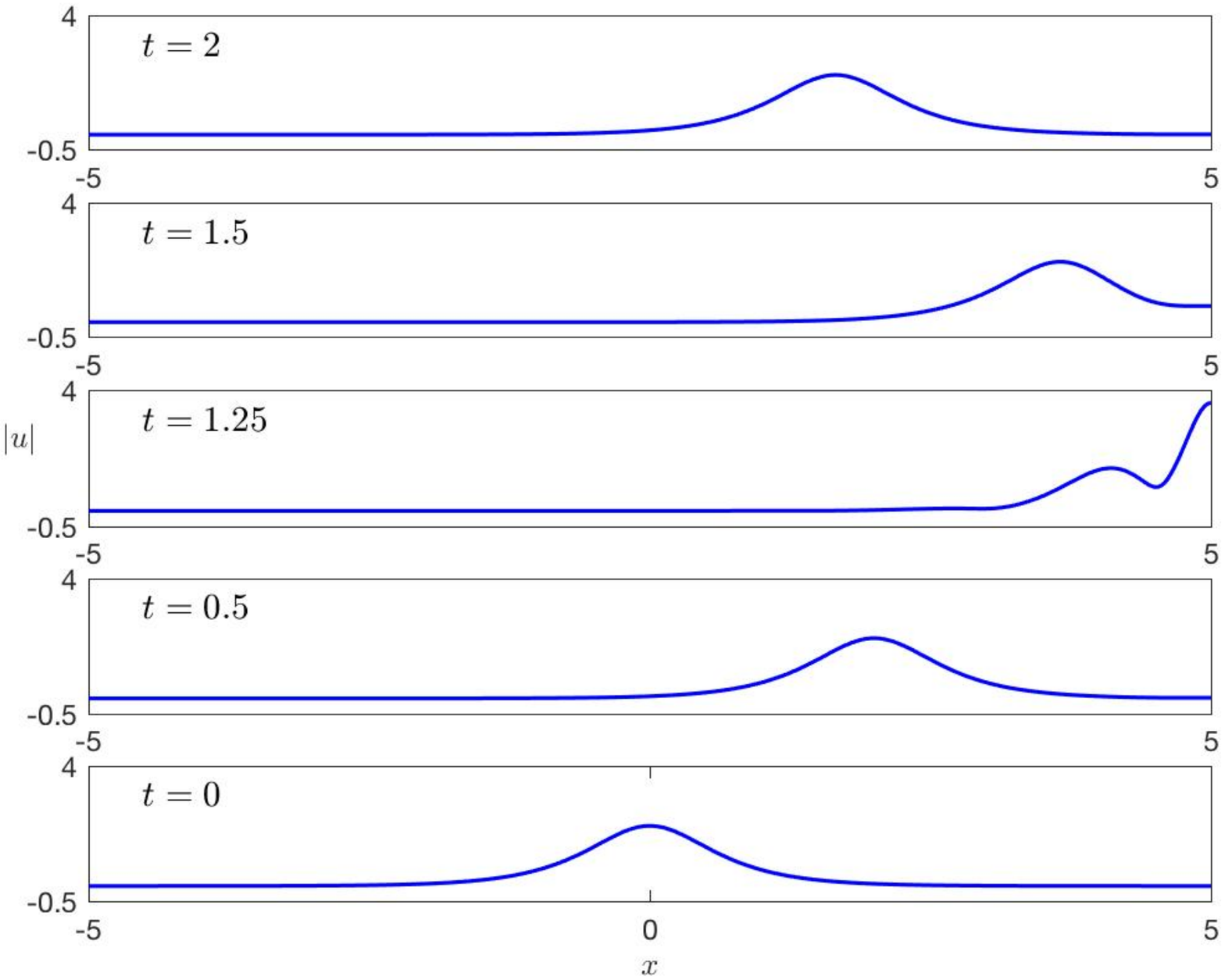}
\caption{Cross section of Fig \ref{refl1a} along the $x$-axis.}
\label{refl1ai}
\end{subfigure}
\begin{subfigure}[b]{0.45\textwidth}
\includegraphics[scale=0.33]{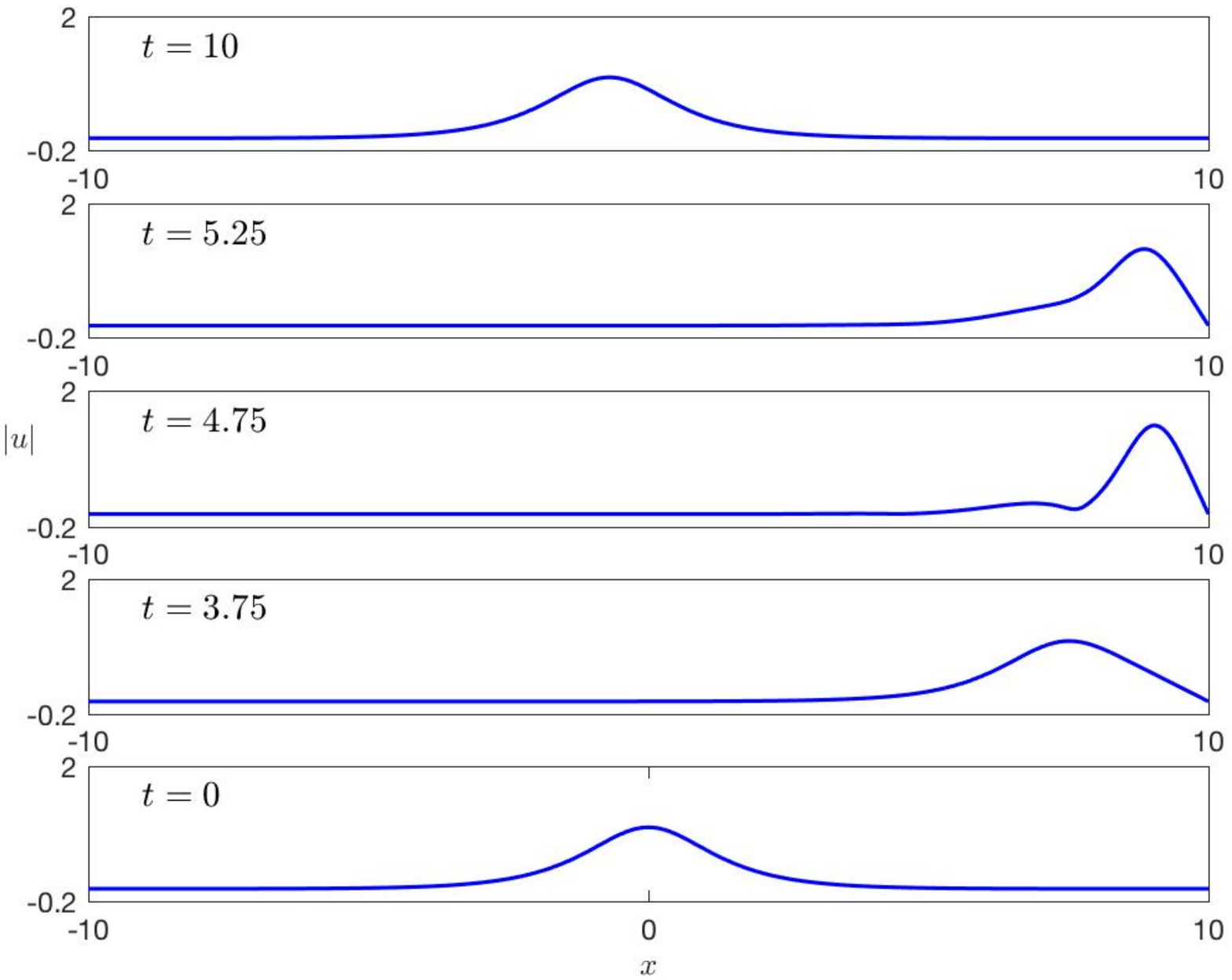}
\caption{Cross section of Fig \ref{refl1b} along the $x$-axis.}
\label{refl1bi}
\end{subfigure} 
\caption{Perfect reflection of bright solitons: cross sections}
\label{refl1i}
\end{figure}

Our second test concerns also the reflection of a bright soliton for the focusing \acs{CNLS} but  in the case where the mixed boundary conditions \eqref{MBC} are used,  c.f. \cite{Biondini1,Biondini2}.  In particular we took $\lambda=2, \ \eta=1$, $\xi=1$ in the domain  $\Omega=[-10,10]\times[-1,1]$ and the zero Dirichlet conditions were applied at $x=-10$ and $x=10$ while homogeneous Neumann boundary conditions in the rest of the boundary of the domain.  We have used an unstructured triangulation consisted of  $74241$ triangles, a time-step $k=2.5\times 10^{-3}$ and quadratic finite elements.  The results of the perfect reflection are presented in Figure \ref{refl1b}, depicting the amplitude of the wave. During the experiment the mass $M$ was conserved with value $3.9999999$ while the energy $E$ was conserved to $1.3333$ up to $T=10$. The reflection in both cases is perfect as it is also can be seen in Figure \ref{refl1i} were the solution is presented along the $x$-axis. The solitons during the interaction with the wall undergoes various changes in their shape and finally regain its original shape and travel in the opposite direction; for more details see \cite{Biondini1,Biondini2}.  

\begin{figure}[ht!]
\centering
\hspace{-2.5cm}
\begin{subfigure}[b]{0.55\textwidth}
\includegraphics[scale=0.35]{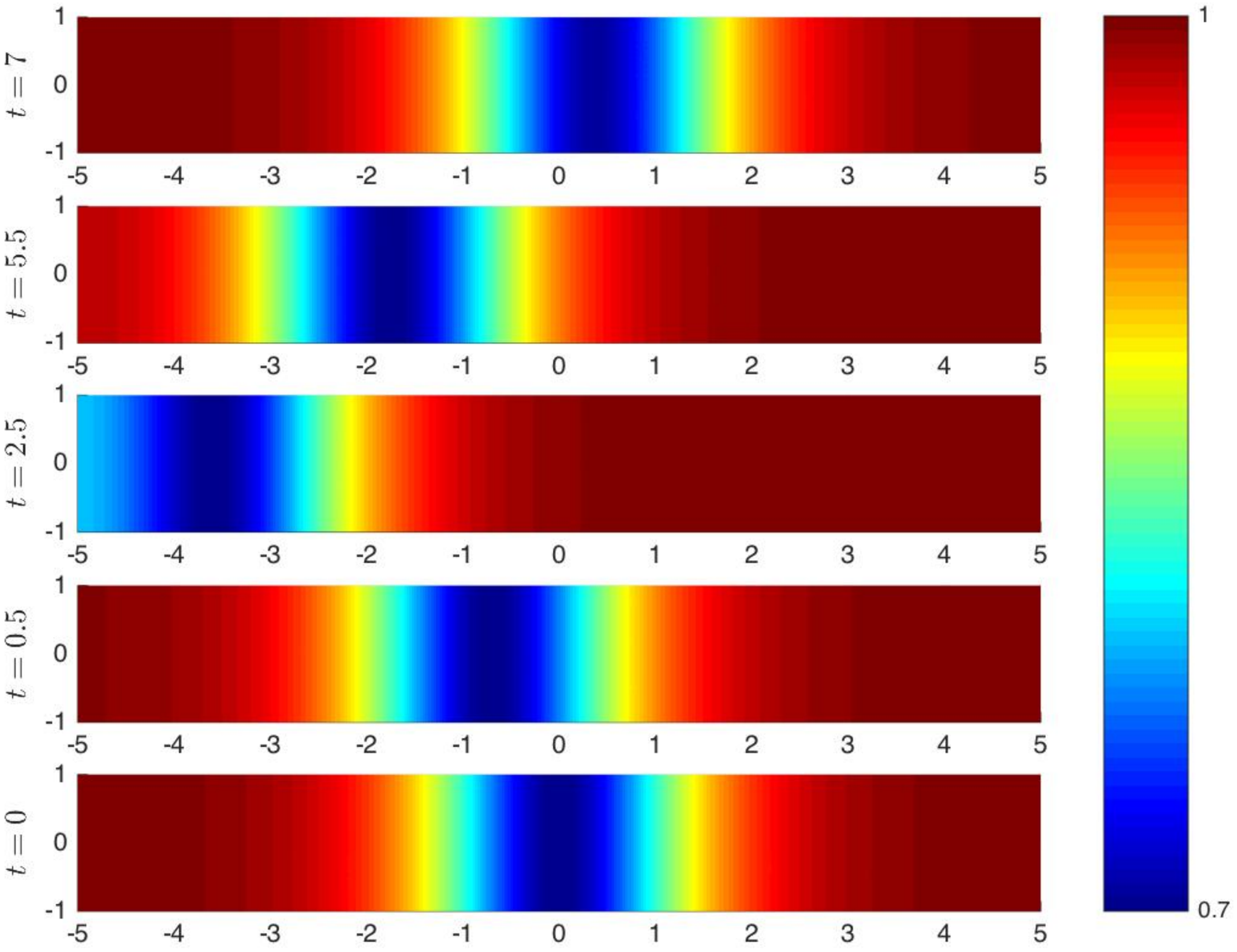}
  \caption{Perfect reflection of a dark soliton.}
  \label{refl2}
\end{subfigure}
\begin{subfigure}[b]{0.45\textwidth}
\includegraphics[scale=0.35]{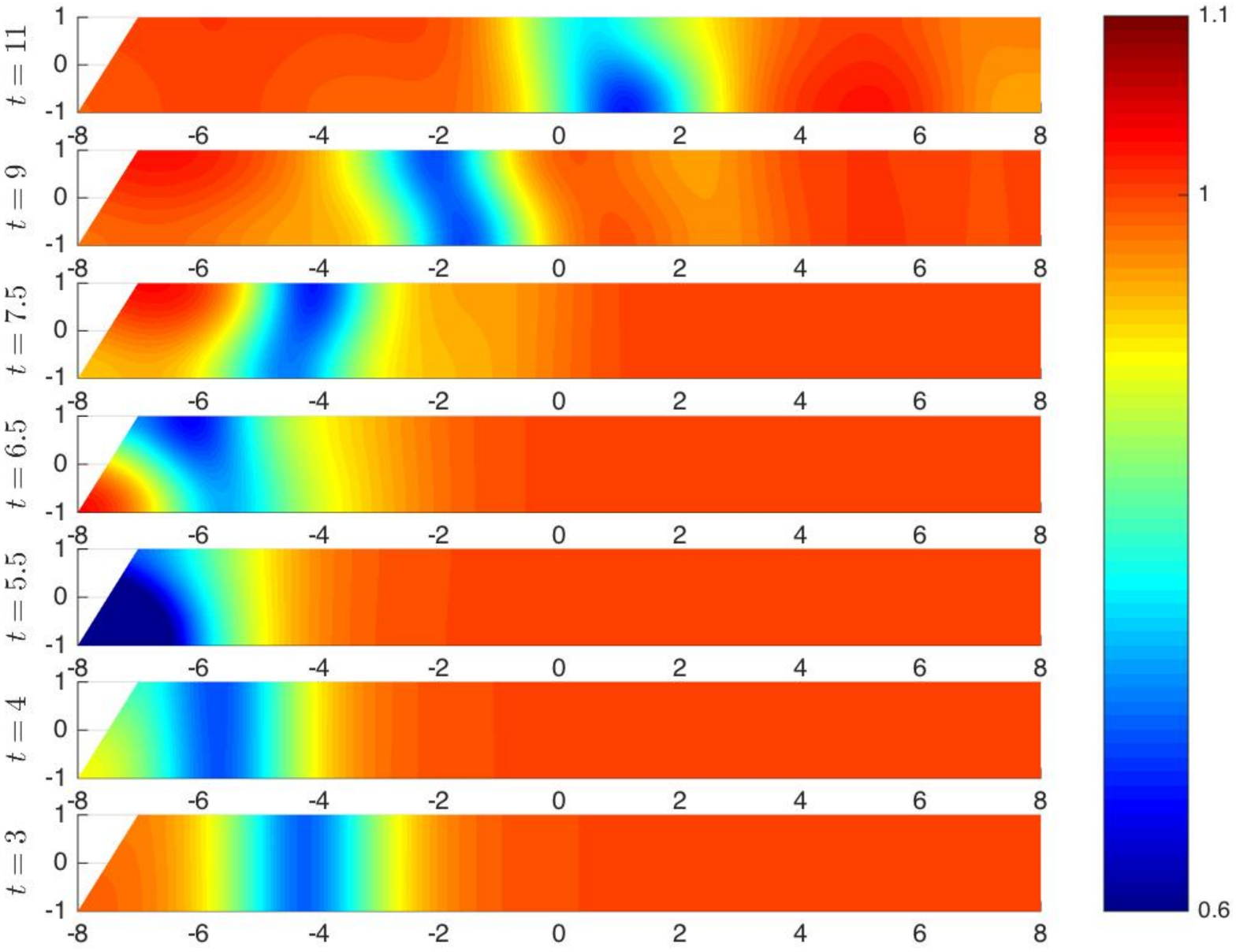}
  \caption{Reflection of a dark soliton by a diagonal wall.}
  \label{refl3}
\end{subfigure} 
\caption{Reflections of dark solitons}
\label{refl4}
\end{figure}

Similarly, the perfect reflection of the dark soliton of the defocusing \acs{CNLS} equation with $\lambda=-2$, $\eta=1$, $\xi=\pi/4$ and with zero Neumann boundary conditions are presented in Figure \ref{refl2}, showing the amplitude of the wave. Since the dark solitons are in general stable waves we used coarser grids than the previous experiment. Specifically, we used an unstructured mesh consisted of $18624$ triangles and time-step $k=5\times 10^{-2}$. Due to the stability properties of the defocusing \acs{CNLS} equation, the mass $M$ was conserved with more digits than in the case of the focusing \acs{CNLS} equation with value $17.1763733098$ while the energy $E$ was conserved to $8.119$ up to $T=15$.

\subsection{Reflection of a dark soliton on a diagonal wall}

The reflection of the dark soliton on a vertical wall presented in the previous section was perfect in the sense that the reflected wave had the same shape with the original soliton. We close this paper with the study of the reflection of the same dark soliton of the defocusing \acs{CNLS} on a diagonal wall. 

Since the \acs{CNLS} equation in a two-dimensional domain is not integrable there is no analytical solution describing such a complicated reflection. The domain that we used here is a trapezoidal domain with vertices $(-8,-1)$, $(-7,1)$, $(8,1)$ and $(8,-1)$ while the triangular grid consisted of $28592$ triangles. In this experiment we use also quadratic finite elements.  The reflection is not perfect and is presented in Figure \ref{refl3}, depicting the amplitude of the wave. As the wave approaches the left side of the boundary, diffraction of the incident wave is being observed. The diffraction of the wave causes the distortion of the soliton while the reflected wave has an oscillatory structure in front and behind of the main pulse. Although it is known that dark solitons exhibit a transverse instability to perturbations with sufficiently long wavelenght, \cite{Kuznetsov, Hoffer}, the reflected wave remained stable and no collapse or any other blow-up phenomenon was observed up to time $T=15$. During this experiment the mass $M$ retained the value (conserving the digits shown) $28.1716833830$ and the energy $E$ was conserved to $13.61$ up to $T=15$. Analogous observations can be made for the focusing \acs{CNLS} equation but are not presented here. 

{\bf Acknowledgements} Dimitrios Mitsotakis would like to thank Profs G. Biondini, S. Flach and B. Ilan for valuable discussions on the properties of the \acs{CNLS} equation. This work was supported by the Victoria University of Wellington Research Establishment Grant (Grand ID 208964). The authors would like to thank the anonymous referees for their valuable comments and suggestions.

\end{document}